\documentclass{amsart}
\usepackage[utf8]{inputenc}
\usepackage[english]{babel}

\usepackage{amssymb,amsthm,amscd,url, hyperref}
\usepackage{amsfonts,amsmath,latexsym} 
\usepackage{amsthm,amscd,mathrsfs,mathtools}
\usepackage{comment}
\usepackage{geometry}
\newcounter{count}
\numberwithin{count}{section}
\newtheorem{theorem}[count]{Theorem}
\newtheorem{lemma}[count]{Lemma}

\newtheorem{corollary}[count]{Corollary}
\newtheorem{definition}{Definition}
\newtheorem{letteredtheorem}{Theorem}

\theoremstyle{definition}
\newtheorem*{theorem*}{Theorem}
\newtheorem*{lemma*}{Lemma}
\newtheorem*{proposition*}{Proposition}
\newtheorem*{question*}{Question}
\newtheorem{remark}[count]{Remark}
\newtheorem*{example*}{Example}

\title[Hutchinson's intervals \& the Laguerre--P\'olya class]{Hutchinson's intervals and 
entire functions from  the Laguerre--P\'olya class}

\author[T.~H.~Nguyen]{Thu Hien Nguyen}
\address{Institute of Mathematics, Leipzig University, Germany, and Department of Mathematics \& Computer Sciences, V.N. Karazin Kharkiv National University,
 Ukraine}
\email{nguyen.hisha@karazin.ua, nguyen.hisha@math.uni-leipzig.de}

\author[A.~Vishnyakova]{Anna Vishnyakova}
\address{Department of Mathematics, Holon Institute of Technology, Israel, and 
Department of Mathematics \& Computer Sciences, V.N. Karazin Kharkiv National University, Ukraine}
\email{anna.vishnyakova@karazin.ua}

\thanks{Aknowledgements: The first author is deeply grateful for the support by 
the Mathematisches Forschungsinstitut Oberwolfach within the program 
Oberwolfach Leibniz Fellows 01.04.2022 -- 30.06.2022.}

\makeatletter
\@namedef{subjclassname@2020}{\textup{2020} Mathematics Subject Classification}
\makeatother



\begin{document}

\begin{abstract}
We find the intervals $[\alpha, \beta (\alpha)]$ such that if a univariate real polynomial 
or entire function $f(z) = a_0 + a_1 z + a_2 z^2 + \cdots $ with positive coefficients satisfy 
the conditions $ \frac{a_{k-1}^2}{a_{k-2}a_{k}} \in [\alpha, \beta(\alpha)]$ for all $k \geq 2,$ 
then $f$ belongs to the Laguerre--P\'olya class. For instance, from J.I.~Hutchinson's theorem, 
one can observe that $f$ belongs to the Laguerre--P\'olya class (has only real zeros) when 
$q_k(f) \in  [4, + \infty).$ We are interested in finding those intervals which 
are not subsets of $[4, + \infty).$
\end{abstract}

\maketitle
  

\section{Introduction}
We study zero localization of real univariate polynomials and entire functions $f(z) = a_0 + a_1 z + a_2 z^2 + \cdots$  with positive coefficients. In 1923, J.I. Hutchinson proved that,  if the inequalities $a_{k-1}^2 
\geq 4 a_{k-2} a_k,$ for all $k \geq 2,$ are valid, then the function $f$ belongs to the Laguerre--P\'olya class.  
In this short note, the chief object is to extend the sufficient conditions for a polynomial or an entire function 
to belong to the Laguerre--P\'olya class obtained by J.I.~Hutchinson, or, more precisely, to find the intervals 
$[\alpha, \beta(\alpha)]$ which are not subsets of $[4, + \infty).$

Let us recall some facts from the theory of entire functions.

\subsection{The Laguerre--P\'olya class}
We begin with the definitions of hyperbolic polynomials, the Laguerre--P\'olya class 
and the Laguerre--P\'olya class of type I.

\begin{definition} A real polynomial $P$ is said to be  \textbf{hyperbolic}, written $P \in \mathcal{HP},$ if all its zeros are real.
\end{definition}

\begin{definition} A real entire function $f$ is said 
to be in the \textbf{Laguerre--P\'olya class}, written $f \in \mathcal{L-P},$ if it can
be expressed in the form
\begin{equation}
\label{lpc}
 f(z) = c z^n e^{-\alpha z^2+\beta z}\prod_{k=1}^\infty
\left(1-\frac {z}{x_k} \right)e^{zx_k^{-1}},
\end{equation}
where $c, \alpha, \beta, x_k \in  \mathbb{R},$ $x_k\ne 0,$  $\alpha \ge 0,$
$n$ is a nonnegative integer and $\sum_{k=1}^\infty x_k^{-2} < \infty.$ 
\end{definition}

\begin{definition} A real entire function $f$ is 
said to be in the \textbf{Laguerre--P\'olya class of type I}, 
written $f \in \mathcal{L-P} I,$ if it can
be expressed in the following form  
\begin{equation} 
\label{lpc1}
 f(z) = c z^n e^{\beta z}\prod_{k=1}^\infty
\left(1+\frac {z}{x_k} \right),
\end{equation}
where $c \in  \mathbb{R},  \beta \geq 0, x_k >0 ,$ 
$n$ is a nonnegative integer,  and $\sum_{k=1}^\infty x_k^{-1} <
\infty.$   
\end{definition}
Note that the product on the right-hand sides in both definitions can be
finite or empty (in the latter case, the product equals 1).

Various important properties and characterizations of the Laguerre--P\'olya class 
and the Laguerre--P\'olya class of type I can be found in works by I.I.~Hirshman and 
D.V.~Widder \cite{HW1955}, B.Ja.~Levin \cite{levin},  G.~P\'olya and G.~Szeg\"o \cite{PSz1998}, 
G.~P\'olya and J.~Schur \cite{PSch1974}, monograph by N.~Obreshkov \cite[Chapter II]{obreshkov} 
and many other works. 
These classes are essential in the theory of entire functions since it appears that the polynomials with 
only real zeros (or only real and nonpositive zeros) converge locally 
uniformly to these and only these functions. The following  prominent theorem provides an 
even stronger result. 

\begin{letteredtheorem}[{E.~Laguerre and G.~P\'{o}lya, see, for example,
\cite[p. 42--46]{HW1955} and \cite[chapter VIII, \S 3]{levin}}]\leavevmode
\label{thmA}
\begin{enumerate}

\item[(i)] Let $(P_n)_{n=1}^{\infty},\  P_n(0)=1, $ be a sequence
of hyperbolic polynomials  which  converges uniformly on the disc 
$|z|\leq A, A > 0.$  Then this sequence converges locally uniformly 
in $\mathbb{C}$ to an entire function from the $\mathcal{L-P}$ class.

\item[(ii)] For any $f \in \mathcal{L-P}$ there exists a sequence of 
hyperbolic polynomials, which converges locally uniformly to $f.$

\item[(iii)] Let $(P_n)_{n=1}^{\infty},\  P_n(0)=1, $ be a sequence
of hyperbolic polynomials having only  negative zeros which  
converges uniformly on the disc $|z| \leq A, A > 0.$  Then this 
sequence converges locally uniformly in $\mathbb{C}$ to an 
entire function  from the class $\mathcal{L-P}I.$ 
 
\item[(iv)] For any $f \in \mathcal{L-P}I$ there is a
sequence of hyperbolic polynomials with only negative 
zeros which converges locally uniformly to $f.$
\end{enumerate}
\end{letteredtheorem}

For a real entire function (not identically zero) of the order less than 
$2$ the property of having only real zeros is  equivalent to belonging 
to the Laguerre--P\'olya class. Similarly, for a real entire function with 
positive coefficients of the order less than $1$ having only real nonpositive 
zeros is  equivalent to belonging to the Laguerre--P\'olya class of type I. 
Strikingly, the situation changes for the functions of order $2$ in the case 
of  the Laguerre--P\'olya class and for the functions of order $1$ in the case 
of  the Laguerre--P\'olya class  of type I.  For instance, 
the entire function $f(x) = e^{-x^2}$ belongs to the $\mathcal{L-P}$ class 
while the entire function $g(x) = e^{x^2}$ does not. 

\subsection{Hutchinson's constant}

The problem of understanding whether a given polynomial or  entire function  has only real zeros 
is considered subtle and complicated. A simply 
verified description of this class, in terms of the 
coefficients of a series, is impossible since it is determined by an infinite number 
of discriminant inequalities. In 1923, J. I. Hutchinson found a simple sufficient 
condition in terms of coefficients for an entire function with positive coefficients 
to have only real zeros, which was a generalization of the results 
by M.~Petrovitch \cite{petrovitch} and G.~Hardy \cite{hardy1904}, or 
\cite[pp. 95 - 100]{hardy1969}.

To formulate the theorem, let us define the second quotients of Taylor coefficients of $f.$
Let  $f(z) = \sum_{k=0}^\infty a_k z^k$  be an entire function with 
real nonzero coefficients, then 
\begin{align}
\label{quotients}
&q_n=q_n(f):=  \frac{a_{n-1}^2}{a_{n-2}a_n}, \quad \forall n\geq 2.
\end{align}
In addition, it follows straightforwardly from this definition that
\begin{align}
\label{q1111}
&a_n  = a_1\Big(\frac{a_1}{a_0} \Big)^{n-1} \frac{1}{q_2^{n-1}q_3^{n-2}
\cdots q_{n-1}^2 q_n}.
\end{align}

\begin{letteredtheorem}[{J.I.~Hutchinson, \cite{hutchinson}}]\label{thmB}
Let $f(z)=
\sum_{k=0}^\infty a_k z^k,$ $ a_k > 0$ for all $k,$ be an entire 
function. Then $q_k(f)\geq 4,$ for all $k\geq 2,$   
if and only if the following two conditions are fulfilled:
\begin{itemize}
\item[(i)] The zeros of $f$ are all real, simple and negative.
\item[(ii)] The zeros of any polynomial $\sum_{k=m}^n a_kz^k,$ $ m < n,$   formed 
by taking any number of consecutive terms of $f ,$ are all real and non-positive.
\end{itemize}
\end{letteredtheorem}
For some extensions of Hutchinson's result see,
for example, the paper by T.~Craven and G.~Csordas, \cite[\S4]{CC1995}.
From Hutchinson's theorem (Theorem~\ref{thmB}) we see that $f$ has 
only real zeros when $q_k(f) \in  [4, + \infty).$

\subsection{Some results related to Hutchinson's constant}

Strikingly, there are many results which are stated in the following style: there 
exists a constant $c > 1$ such that if a  polynomial or an entire function $f$ 
with nonzero coefficients  satisfies the conditions $|q_k(f)| \geq c$ for all $k,$ 
then we can formulate a statement about the localization of the zeros of $f.$ 
For example, in \cite{CC1995} the authors obtained an analogue of the Hutchinson's 
theorem for polynomials decomposed in the Pochhammer basis. In \cite{handelman}, it was proved that, if  for some constant $c>0$ a  
polynomial $P$ with positive coefficients satisfies the conditions $q_k(P) > c$ for all $k,$ then all the zeros of $P$ lie in a special sector depending on $c.$   In \cite{KV2008}, the 
smallest possible constant $c>0$ was found such that  if a polynomial $P$ with positive
coefficients satisfies the conditions $q_k(P) > c$ for all $k,$ then $P$ is stable (all the zeros 
of $P$ lie in the left half-plane). In \cite{KarVi2017}, the smallest possible constant $c>0$ was 
found such that  if a polynomial $P$ with positive coefficients satisfies the conditions $q_k(P) > c$ 
for all $k,$ then $P$ is a sign-independently hyperbolic polynomial.  In  \cite{BNV},
the smallest possible constant $c>0$ was found such that  if a polynomial $P$ with 
complex coefficients satisfies the conditions $|q_k(P)| > c$ 
for all $k,$ then $P$ has only simple zeros.

The following special function 

$$g_a(z) =\sum _{k=0}^{\infty} z^k a^{-k^2}, \ a > 1,$$ 
which is called the \textbf{partial theta function}, plays a significant role in the mentioned circle of problems. Strikingly, $ q_k(g_a) = a^2 $ for all $ k \geq 2.$ One of the interesting questions is, for which values of $a$  this function belongs to the Laguerre--P\'olya class. The paper \cite{KLV2003} by  O.M.~Katkova, T.~Lobova-Eisner, and A.M.~Vishnyakova gives 
an exhaustive answer to this question. In particular, it is proved that there exists a constant $ q_\infty\approx 3{.}23363666$
such that $g_a \in \mathcal{L-P}$   if and only if $a^2 \geq q_\infty.$ Moreover, the authors studied analogous questions 
for the Taylor sections of the function $g_a.$  For more details on the partial theta function, 
see a series of works by V.P.~Kostov dedicated to its various properties \cite{kostov2013a, kostov2013b,    
kostov2015, kostov2016b}, his joint work with B.~Shapiro \cite{KSh2013}, and a fascinating historical review by S.O.~Warnaar \cite{warnaar2018}.

It is easy to show that, if the estimation of $q_k(f)$ only from below is
given then the constant  $4$  in $q_k (f) \geq 4$ is the smallest possible to conclude
that $f  \in \mathcal{L-P}$ (that is, Theorem~\ref{thmB} remains valid when omitting (ii)). However, if we only have the estimation of $q_k$ from below and  require monotonicity, then the constant $4$ in the condition $q_k \geq 4$ can be reduced to conclude that $f \in \mathcal{L-P}.$  As an example, in \cite{NgV2018}, it was proved that if the entire functions have the decreasing $q_k$  such that $\lim_{n \to \infty} q_k = c \geq q_\infty,$  then the function belongs to the Laguerre--P\'olya class.

In this work, we show that if the estimations on $q_k(f)$ from below and from  above are given, then the constant $ 4$ can be decreased. We would like to investigate such problems where assumption 
$q_k(f) \geq c$ for all $k$ is changed by $q_k(f) \in [\alpha, \beta]$ for all $k$ for 
some given segment $[\alpha, \beta].$ As far as we know, the first result of such kind 
was obtained in  \cite{KLV2003} where the following theorem was proved.

\begin{letteredtheorem}[{O.M.~Katkova, T.~Lobova, and A.M.~Vishnyakova, \cite{KLV2003}}]\label{thmC}
Let $f(z)=\sum_{k=0}^\infty a_k z^k,$ $a_k>0, $ be an entire function and $\alpha \in [3{.}43 ; 4].$
Then $q_k(f)  \in \left[\alpha , \frac {0{.}95}{2\sqrt \alpha -\alpha} \right]$ for all $k \geq 2$ implies 
$f \in \mathcal{L-P}.$
\end{letteredtheorem}

\section{Hutchinson's intervals}

We present our main result.

\begin{theorem}
\label{thm1}
Let $P(x) = \sum_{k=0}^n a_k x^k,$ $a_k > 0,$  be a polynomial, and $n \geq 4.$ Suppose that 
there exists $\alpha,  1 + \sqrt{5} \leq \alpha < 4,$ such that
 $q_k(P) \in \left[\alpha, \frac{8}{\alpha(4 - \alpha)} \right]$ for all $k =2, 3, \ldots, n.$  
Then $P \in \mathcal{HP}.$
\end{theorem}

The following statement is a simple corollary of the above result.

\begin{corollary}
\label{Cor1}
Let $f(x) = \sum_{k=0}^\infty a_k x^k,$ $a_k > 0,$  be an entire function. Suppose that 
there exists $\alpha,  1 + \sqrt{5} \leq \alpha < 4,$ such that
 $q_k(f) \in \left[\alpha, \frac{8}{\alpha(4 - \alpha)}\right]$ for all $k =2, 3, \ldots$  
Then $f \in \mathcal{L-P}.$
\end{corollary}

\begin{remark}
As it follows from the result about the partial theta-function by \cite{KLV2003},
the constant  $\alpha$ in the statement of Theorem \ref{thm1} can
not be less than $ q_\infty \approx 3{.}23363666.$ We observe that
$1 + \sqrt{5} \approx  3.23606797.$

\end{remark}

\subsection{Proof of Theorem~\ref{thm1}}

For a polynomial $P(x) = \sum_{k=0}^n a_k x^k$  with  positive coefficients, without loss of generality, we can assume that $a_0=a_1=1,$ since we can 
consider a polynomial $T(x) = a_0^{-1} P (a_0 a_1^{-1}x) $  instead of 
$P(x),$ due to the fact that such rescaling of $P$ preserves its property of 
having real zeros and preserves the second quotients:  $q_k(T) =q_k(P)$ 
for all $k.$ For the sake of brevity, we further use notation  $q_k$ instead 
of  $q_k(P).$ Thereafter, we consider a polynomial  
\begin{align}
\label{phi}
Q(x) = T(-x) = 1 - x + \sum_{k=2}^n  \frac{ (-1)^k x^k}
{q_2^{k-1} q_3^{k-2} \cdots q_{k-1}^2 q_k}
\end{align}
instead of $P$ (see (\ref{q1111}) for the formulas for coefficients).

Our proof is based on the following lemma.

\begin{lemma}
\label{thm1:lm1}
Let $[\alpha, \beta],$ $0 < \alpha < \beta,$ be a given segment.
Then the following two statements are equivalent:
\begin{itemize}
\item[(a)] For every polynomial
\begin{align*}
S_{q_2, q_3, q_4}(x) = 1 - x + \frac{x^2}{q_2} - \frac{x^3}{q_2^2 q_3} + \frac{x^4}{q_2^3 q_3^2 q_4}
\end{align*}
such that $q_j \in [\alpha, \beta]$ for all $j =2, 3, 4,$ there 
exists a point $x_0 \in (1, \alpha)$ such that $S_{q_2, q_3, q_4}(x_0) \leq 0.$ 
\item[(b)] The following inequalities are valid: $\alpha \geq  1 + \sqrt{5},$ and, if  $\alpha < 4$ then 
$\beta \leq \frac{8}{\alpha(4 - \alpha)}.$
\end{itemize}
\end{lemma}

\begin{proof}

Suppose that the Statement~(a) of Lemma~\ref{thm1:lm1}  is valid. Then for the polynomial 
$S_{\alpha, \alpha, \alpha}(x)$ $= 1 - x + \frac{x^2}{\alpha} - 
\frac{x^3}{\alpha^3} + \frac{x^4}{\alpha^6},$ 
 there exists a point $x_0 \in (1, \alpha)$ such that $S_{\alpha, \alpha, \alpha}(x_0) \leq 0.$ It is easy to
 check the following identity
\begin{equation}
\label{an1}
S_{\alpha, \alpha, \alpha}(x) =\left(\frac{x^2}{\alpha^3}  - \frac{x}{2}+1\right)^2 - 
\left(\frac{1}{4} + \frac{2}{\alpha^3} 
- \frac{1}{\alpha}\right) x^2.
\end{equation}
If $\left(\frac{1}{4} + \frac{2}{\alpha^3} 
- \frac{1}{\alpha}\right) <0,$ then for every $x \in \mathbb{R}$ we have $S_{\alpha, \alpha, \alpha}(x) > 0.$
Thus, $$\left(\frac{1}{4} + \frac{2}{\alpha^3} 
- \frac{1}{\alpha}\right) = \frac{1}{4 \alpha^3}(\alpha -1 + \sqrt{5})(\alpha -2)
(\alpha -1 - \sqrt{5})  \geq 0,$$ whence $ \alpha \in (0, 2] \cup [1 + \sqrt{5}, \infty).$

First we consider the case $ \alpha \in (0, 2].$ By (\ref{an1}), we have
\begin{align*}
S_{\alpha, \alpha, \alpha}(x) &=\left(\frac{x^2}{\alpha^3} -  \left(\frac{1}{2}   + 
\sqrt{\frac{1}{4} + \frac{2}{\alpha^3} 
- \frac{1}{\alpha}} \right) x +1\right) \\
&\quad \times \left(\frac{x^2}{\alpha^3} -  \left(\frac{1}{2}   - 
\sqrt{\frac{1}{4} + \frac{2}{\alpha^3} 
- \frac{1}{\alpha}} \right) x +1\right).
\end{align*}
We have two quadratic polynomials in brackets with the following discriminants:
$$ D_{\pm} = \left(\frac{1}{2}   \pm \sqrt{\frac{1}{4} + \frac{2}{\alpha^3} 
- \frac{1}{\alpha}} \right)^2 - \frac{4}{\alpha^3}.$$
If $D_{+} <0$ and $D_{-} <0,$ then for every $x \in \mathbb{R}$ we have 
$S_{\alpha, \alpha, \alpha}(x) > 0.$
Thus, at least one of these two discriminants is nonnegative, whence $D_{+} \geq 0,$
and we obtain
\begin{align}
\label{ineq1}
 \sqrt{\frac{1}{4} + \frac{2}{\alpha^3} 
- \frac{1}{\alpha}} \geq \frac{4 - \alpha^{3/2}}{2\alpha^{3/2}}.
\end{align}
We observe that $4 - \alpha^{3/2} \geq 0$ for $\alpha \in (0, 2].$ Thus, the inequality (\ref{ineq1}) implies

\begin{align*}
 \psi(\alpha):= -\alpha^2 + 2 \alpha^{3/2} - 2  \geq 0.  
 \end{align*}
The derivative of $\psi(\alpha)$
has a unique positive root  $\alpha_0 = \frac{9}{4},$
and the maximal value of $\psi$ for $\alpha >0$ is
$\psi\left(\frac{9}{4}\right) = - \frac{5}{16} < 0.$ Thus, if 
the statement~(a) of Lemma~\ref{thm1:lm1} is valid then $\alpha \geq  1 + \sqrt{5}.$

Further, we assume that $\alpha \geq  1 + \sqrt{5},$
and suppose that the statement~(a) of Lemma~\ref{thm1:lm1} is valid. Let 
$S_{q_2, q_3, q_4}(x)$ be an arbitrary polynomial
such that $q_j \in [\alpha, \beta]$ for all $j =2, 3, 4.$ 
We want to investigate whether there exists $x_0 \in (1, \alpha)$ 
such that $S_{q_2, q_3, q_4}(x_0) \leq 0.$ 
We observe that for all $x>0$
\begin{align*} 
S_{q_2, q_3, q_4}(x) \leq S_{q_2, q_3, \alpha}(x) := 1 - x + 
\frac{x^2}{q_2} - \frac{x^3}{q_2^2 q_3} + \frac{x^4}{q_2^3 q_3^2 \alpha}. 
\end{align*} 
Thus, for every polynomial $S_{q_2, q_3, q_4}(x)  $
with $q_j \in [\alpha, \beta]$ for all $j =2, 3, 4,$ there 
exists a point $x_0 \in (1, \alpha)$ such that $S_{q_2, q_3, q_4}(x_0) \leq 0$ 
if and only if for every polynomial $S_{q_2, q_3, \alpha}(x)$
with $q_j \in [\alpha, \beta]$ for all $j =2, 3,$  there 
exists a point  $x_0 \in (1, \alpha)$ such that $S_{q_2, q_3, \alpha}(x_0) \leq 0.$

Next, we compute the derivative of $S_{q_2, q_3, \alpha}$ with respect to $q_3.$ We get
\begin{align*}
  \frac{\partial}{\partial q_3}S_{q_2, q_3, \alpha} (x) = \frac{x^3}{q_2^2 q_3^2} - 
\frac{2x^4}{q_2^3 q_3^3 \alpha}.
\end{align*}
We observe that $\frac{x^3}{q_2^2 q_3^2} - 
\frac{2x^4}{q_2^3 q_3^3 \alpha} >0 \Leftrightarrow x< \frac{q_2 q_3 \alpha}{2},$  
so for all $x \in (1, \alpha)$ we get that $S_{q_2, q_3, \alpha}$ is increasing in $q_3.$ 
Whence we have
\begin{align*}
S_{q_2, q_3, \alpha}(x) \leq S_{q_2, \beta, \alpha}(x) = 1 - x + \frac{x^2}{q_2} - 
\frac{x^3}{q_2^2 \beta} + \frac{x^4}{q_2^3 \beta^2 \alpha}. 
\end{align*}
Analogously, we consider the derivative of $S_{q_2, \beta, \alpha}(x)$ with respect to $q_2$ 
to understand the monotonicity and we get
\begin{align*}
  \frac{\partial}{\partial q_2}S_{q_2, \beta, \alpha} (x) = -\frac{x^2}{q_2^2} + 
\frac{2x^3}{q_2^3 \beta} - \frac{3x^4}{q_2^4 \beta^2 \alpha}.
\end{align*}
We show that $\frac{\partial}{\partial q_2}S_{q_2, \beta, \alpha} (x) < 0$ for $x \in (1, \alpha),$
or, equivalently,
\begin{align}
\label{ineq2}
3 x^2 -2 q_2 \alpha \beta  x + q_2^2  \alpha \beta^2>0.
\end{align}
Under our assumption that $\alpha \geq 1 + \sqrt{5},$ we compute the discriminant of the lefthand side of (\ref{ineq2})and observe that
$\frac{D}{4} = q_2^2\beta^2 \alpha (\alpha -3) > 0,$
so the quadratic expression has two
positive roots $$x_{\pm} = \frac{q_2 \alpha \beta \pm q_2 \beta \sqrt{\alpha 
(\alpha -3)}}{3}.$$ To prove (\ref{ineq2}),
it is sufficient to check that $\alpha < x_{-},$  or  $q_2 \beta \sqrt{\alpha 
(\alpha -3)} < q_2 \alpha \beta -3 \alpha.$ The last inequality is equivalent 
to  $q_2^2 \beta^2 +3\alpha - 2 q_2 \alpha \beta = q_2\beta (q_2\beta -2\alpha)
+3\alpha >0,$ which holds under our assumptions since $q_2 \geq \alpha,$ and $\beta > \alpha >2.$
Thus, we have proved that for all $x \in (1, \alpha)$:
$$ S_{q_2, \beta, \alpha} (x) \leq S_{\alpha, \beta, \alpha} (x) 
= 1 - x + \frac{x^2}{\alpha} - \frac{x^3}{\alpha^2\beta} + \frac{x^4}{\alpha^4\beta^2}.$$
Consequently, for every polynomial $S_{q_2, q_3, q_4}(x)  $
such that $q_j \in [\alpha, \beta]$ for all $j =2, 3, 4,$ there 
exists a point $x_0 \in (1, \alpha)$ such that $S_{q_2, q_3, q_4}(x_0) \leq 0$ 
if and only if for the polynomial $S_{\alpha, \beta, \alpha}(x) $
 there exists a point $x_0 \in (1, \alpha)$ such that $S_{\alpha, \beta, \alpha}(x_0) \leq 0.$

Now we consider the polynomial $S_{\alpha, \beta, \alpha} (x)$  for $x \in (1, \alpha).$
Set $x =: \alpha\sqrt{\beta} y.$ Since $x \in (1, \alpha),$
we have $y \in (\frac{1}{\alpha\sqrt{\beta}}, \frac{1}{\sqrt{\beta}}).$
Hence, after change of variables, we get a self-reciprocal polynomial
\begin{align*}
\widetilde{P}(y) := S_{\alpha, \beta, \alpha} ( \alpha\sqrt{\beta} y) &= 1 -  
\alpha\sqrt{\beta} y + \alpha \beta y^2 - \alpha\sqrt{\beta} y^3 + y^4\\
&=y^2 \bigg( (y^{-2} + y^2 ) - \alpha\sqrt{\beta} (y^{-1} + y) + \alpha \beta \bigg).
\end{align*}
Set $w := y^{-1} + y.$ 
We want to investigate whether there exists a point $w_0 \in (\sqrt{\beta} + \frac{1}{\sqrt{\beta}}, 
\alpha\sqrt{\beta} + \frac{1}{\alpha\sqrt{\beta}})$ such that
\begin{align*}
\widetilde{\widetilde{P}}(w_0) = w_0^2 - \alpha\sqrt{\beta}w_0 + \alpha \beta - 2 \leq 0.
\end{align*}
We consider the vertex of the parabola
$w_{v} = \frac{\alpha\sqrt{\beta}}{2},$
and check if it lies in $\left(\sqrt{\beta} + \frac{1}{\sqrt{\beta}}, 
\alpha\sqrt{\beta} + \frac{1}{\alpha\sqrt{\beta}} \right).$ Obviously,  	
$w_{v} = \frac{\alpha\sqrt{\beta}}{2} < \alpha\sqrt{\beta} + 
\frac{1}{\alpha\sqrt{\beta}}.$  We show that the following inequality is fulfilled 
$ \frac{\alpha\sqrt{\beta}}{2} > \sqrt{\beta} + \frac{1}{\sqrt{\beta}},$
or, equivalently, $\beta > \frac{2}{\alpha-2}.$
It is sufficient to prove that
$\alpha > \frac{2}{\alpha-2},$
and it is equivalent to
$\alpha^2 - 2\alpha - 2 > 0,$
which is fulfilled under our assumption  $\alpha \geq 1 + \sqrt{5}.$
Since $w_{v} \in \left(\sqrt{\beta} + \frac{1}{\sqrt{\beta}}, 
\alpha\sqrt{\beta} + \frac{1}{\alpha\sqrt{\beta}} \right),$  there exists 
$w_0 \in \left(\sqrt{\beta} + \frac{1}{\sqrt{\beta}}, 
\alpha\sqrt{\beta} + \frac{1}{\alpha\sqrt{\beta}}\right)$ such that
$\widetilde{\widetilde{P}}(w_0)  \leq 0$
if and only if  the discriminant of this quadratic function is non-negative:
\begin{align}
\label{suffcond}
D = \alpha^2 \beta - 4 \alpha \beta + 8 = \beta \alpha ( \alpha -4) +8\geq 0.
\end{align}
The inequality above is equivalent to the following statement: if  $\alpha < 4$ then 
$\beta \leq \frac{8}{\alpha(4 - \alpha)}.$
\end{proof}

\begin{remark}
Lemma~\ref {thm1:lm1} is an analog of
Theorem 1.5 from \cite{NgV2020}. This theorem states that
if $f(x) = 1 + x + \sum_{k=2}^\infty a_k x^k$ is an entire function with 
positive coefficients, and  $3  \leq  q_2(f) < 4, $  $q_4(f) \geq 3$ and
$2 \leq q_3(f) \leq \frac{8}{d(4-d)},$
where  $d = \min(q_2(f), q_4(f)),$ then  there exists $x_0 \in [-q_2(f), 0]$ 
such that $f(z_0) \leq 0.$  
\end{remark}


Now we can prove Theorem~\ref{thm1}. 
Let $Q(x) =  1 - x + \sum_{k=2}^n  \frac{ (-1)^k x^k}
{q_2^{k-1} q_3^{k-2} \cdots q_{k-1}^2 q_k}$  be a polynomial,  $n \geq 4,$ and 
there exists $\alpha \in [1 + \sqrt{5}, 4)$ such that
 $q_k \in \left[\alpha, \frac{8}{\alpha(4 - \alpha)} \right]$ for all $k =2, 3, \ldots, n.$ 
Let us fix an arbitrary $j$ such that $1 \leq j \leq \lfloor\frac{n}{2}\rfloor,$  and suppose that 
$x \in  (q_2 q_3 \cdots q_j, q_2 q_3 \cdots q_j q_{j+1})$ (for $j=1$ we assume that 
$1 < x < q_2$). Then we observe that
\begin{align}
\label{monot1}
1 < x < \frac{x^2}{q_2} < \frac{x^3}{q_2^2 q_3} < \cdots < \frac{x^j}{q_2^{j-1} q_3^{j-2} \cdots q_j},
\end{align}
and
\begin{align}
\label{monot2}
 \frac{x^j}{q_2^{j-1} q_3^{j-2} \cdots q_j} &> \frac{x^{j+1}}{q_2^{j} q_3^{j-1} 
\cdots q_j^2 q_{j+1}} > \\ &  \nonumber   \frac{x^{j+2}}{q_2^{j+1} q_3^{j} \cdots q_j^3 
q_{j+1}^2 q_{j+2}} > \cdots  > \frac{x^{n}}{q_2^{n-1} q_3^{n-2} \cdots  
q_{n-1}^2 q_{n}}.
\end{align}
We have the following representation:
\begin{align}
\label{dissec1}
  (-1)^{j-1} Q(x) & =  \left((-1)^{j-1} - (-1)^{j-1} x + \sum_{k=2}^{j-2}  \frac{ (-1)^{k+j-1} x^k}
{q_2^{k-1} q_3^{k-2} \cdots q_{k-1}^2 q_k}\right) 
 \\   \nonumber  & \quad +  \sum_{k=j-1}^{j+3}  \frac{ (-1)^{k+j -1} x^k}
{q_2^{k-1} q_3^{k-2} \cdots q_{k-1}^2 q_k} + \left(\sum_{k=j+4}^n  \frac{ (-1)^{k+j-1} x^k}
{q_2^{k-1} q_3^{k-2} \cdots q_{k-1}^2 q_k}\right) 
 \\   \nonumber  & := \Sigma_{1,j}(x) + g_j(x) + \Sigma_{2,j}(x).  
\end{align}
We note that, for some $j$ the sum $\Sigma_{2, j}(x)$ can be empty (and equal to zero), but for $n\geq 5$ we have $j +3  \leq \lfloor\frac{n}{2}\rfloor +3 \leq n,$ 
so all $5$ summands in $g_j(x)$ are nonzero. We later consider the case $n=4.$

For  $x \in  (q_2 q_3 \cdots q_j, q_2 q_3 \cdots q_j q_{j+1}),$ 
we observe that the terms in $\Sigma_{1, j}(x)$ are alternating in sign and and their moduli are increasing, 
while the summands in $\Sigma_{2, j}(x)$ are alternating in sign and their moduli are decreasing. Hence, 
$\Sigma_{1, j}(x) < 0$ and $\Sigma_{2, j}(x) < 0$ for  all  $x \in  (q_2 q_3 \cdots q_j, q_2 q_3 
\cdots q_j q_{j+1}),$ whence we get
\begin{align}
\label{est1}
(-1)^{j-1} Q(x)  <  g_j(x) \quad \forall x \in   (q_2 q_3 \cdots q_j, q_2 q_3 
\cdots q_j q_{j+1}).
\end{align}
We have 
\begin{align*}
g_j(x)  &= \frac{ x^{j-1}}{q_2^{j-2} q_3^{j-3} \cdots q_{j-2}^2 q_{j -1}} \bigg(1 - \frac{x}{q_2 q_3 \cdots q_{j-1}q_j}  
+ \frac{x^2}{q_2^2 q_3^2 \cdots q_{j-1}^2q_j^2 q_{j+1}}  \\
 &\quad -  \frac{x^3}{q_2^3 q_3^3 \cdots q_{j-1}^3 q_j^3 q_{j+1}^2 q_{j+2}} 
+\frac{x^4}{q_2^4 q_3^4 \cdots q_{j-1}^4 q_j^4 q_{j+1}^3 q_{j+2}^2 q_{j+3}}  \bigg).  
\end{align*}
Set $y =\frac{x}{q_2 q_3 \cdots q_{j-1}q_j},$ and for $x \in  (q_2 q_3 \cdots q_j, q_2 q_3 
\cdots q_j q_{j+1})$ we have $y \in (1, q_{j+1}).$ Therefore, we obtain
\begin{align*}
h_j(y) :=   g_j(q_2 q_3 \cdots q_{j-1}q_j y)  = q_2 q_3^{2} \cdots q_{j -1}^{j-2}
q_j^{j-1} y^{j-1} \\
\times \bigg(1 - y + \frac{y^2}{q_{j+1}}  -   
 \frac{y^3}{q_{j+1}^2q_{j+2}}
+ \frac{y^4}{q_{j+1}^3q_{j+2}^2 q_{j+3}}\bigg).  
\end{align*}
Since $q_{j+1}, q_{j+2}, q_{j+3} \in \left[\alpha, \frac{8}{\alpha(4 - \alpha)}\right] $ for some $\alpha,$ where 
 $1 + \sqrt{5} \leq \alpha < 4,$ the polynomial in brackets satisfies the assumptions of 
Lemma~\ref {thm1:lm1}. Thus, there exists $y_j \in  (1, q_{j+1})$ such that
$h_j(y_j) \leq 0.$ Hence, there exists $x_j = q_2 q_3 \cdots q_{j-1}q_j y_j \in 
(q_2 q_3 \cdots q_j, q_2 q_3 
\cdots q_j q_{j+1}),$ such that $ g_j(x_j)\leq 0.$ Taking into account (\ref{est1}),
 for every $n \geq 5$ we obtain
\begin{equation}
\label{signs1}
\forall j, \   1 \leq j \leq \lfloor \frac{n}{2}\rfloor,\quad \exists 
x_j  \in (q_2 q_3 \cdots q_j, q_2 q_3 
\cdots q_j q_{j+1})  : \quad (-1)^{j-1} Q(x_j) <0.
\end{equation}

The only problem for $n=4$ is when $j=2.$ In the latter case,
we have
\begin{align*}
 - Q(x) = -1 + \left( x - \frac{x^2}{q_2} + \frac{x^3}{q_2^2 q_3} - 
\frac{x^4}{q_2^3 q_3^2 q_4}\right) = -1 + x\left( 1 - \frac{x}{q_2} + \frac{x^2}{q_2^2 q_3} - 
\frac{x^3}{q_2^3 q_3^2 q_4}\right).
\end{align*}
We highlight that the polynomial in brackets is of degree $3,$ however, we can estimate it from above with a polynomial of degree $4$ as follows:
\begin{align*}
 - Q(x) < x\left( 1 - \frac{x}{q_2} + \frac{x^2}{q_2^2 q_3} - 
\frac{x^3}{q_2^3 q_3^2 q_4}\right) < x\left( 1 - \frac{x}{q_2} + \frac{x^2}{q_2^2 q_3} - 
\frac{x^3}{q_2^3 q_3^2 q_4} + \frac{x^4}{q_2^4 q_3^3 q_4^2 q_4}\right).
\end{align*}
Therefore, we can further reason in the same way as before. Thus, we conclude that (\ref{signs1}) is valid for all $n \geq 4.$

Now let us fix an arbitrary $j,  \lfloor\frac{n}{2}\rfloor +1 \leq j \leq n-1,$  and suppose that 
$x \in$ $(q_2 q_3 \cdots q_j, q_2 q_3 \cdots q_j q_{j+1}).$ Then both (\ref{monot1}) and 
(\ref{monot2})  are valid, and we have the following representation:
\begin{align}
\label{dissec2}
  (-1)^{j-1} Q(x) &=  \left((-1)^{j-1} - (-1)^{j-1} x + \sum_{k=2}^{j-4}  \frac{ (-1)^{k+j-1} x^k}
{q_2^{k-1} q_3^{k-2} \cdots q_{k-1}^2 q_k}\right) 
 \\   \nonumber  &\quad +  \sum_{k=j-3}^{j+1}  \frac{ (-1)^{k+j -1} x^k}
{q_2^{k-1} q_3^{k-2} \cdots q_{k-1}^2 q_k} + \left(\sum_{k=j+2}^n  \frac{ (-1)^{k+j -1} x^k}
{q_2^{k-1} q_3^{k-2} \cdots q_{k-1}^2 q_k}\right) 
 \\   \nonumber  &
:= \Sigma_{1,j}(x) + g_j(x) + \Sigma_{2,j}(x).  
\end{align}
We note that for some $j,$ the sum $\Sigma_{1,j}(x)$ can be empty (and
equal to zero), although for $n\geq 4,$ we have $j -3  \geq \lfloor\frac{n}{2}\rfloor -2 \geq 0,$ 
so all $5$ summands in $g_j(x)$ are nonzero. 

For  $x \in  (q_2 q_3 \cdots q_j, q_2 q_3 \cdots q_j q_{j+1}),$ 
we obeserve that the terms in $\Sigma_{1,j}(x)$ are alternating in sign and and their moduli are increasing, 
while the summands in $\Sigma_{2,j}(x)$ are alternating in sign and their moduli are decreasing. Hence, 
$\Sigma_{1,j}(x) < 0$ and $\Sigma_{2,j}(x) < 0$ for  all  $x \in  (q_2 q_3 \cdots q_j, q_2 q_3 
\cdots q_j q_{j+1}),$ whence we get
\begin{align}
\label{est2}
(-1)^{j-1} Q(x)  <  g_j(x) \quad \forall x \in   (q_2 q_3 \cdots q_j, q_2 q_3 
\cdots q_j q_{j+1}).
\end{align}
We have 

\begin{align*}
g_j(x)  = \frac{ x^{j+1}} {q_2^{j} q_3^{j-1} \cdots q_{j}^2 q_{j +1}} &\times \bigg(1 - \frac{q_2 q_3 \cdots q_{j-1}q_j q_{j+1}}{x}  + \frac{q_2^2 q_3^2 \cdots q_{j-1}^2q_j^2 q_{j+1}}{x^2}  
 \\ & \quad -  \frac{q_2^3 q_3^3 \cdots q_{j-2}^3q_{j-1}^3 q_j^2 q_{j+1}}{x^3} 
+\frac{q_2^4 q_3^4 \cdots q_{j-3}^4 q_{j-2}^4q_{j-1}^3 q_j^2 q_{j+1}}{x^4}  \bigg).  
\end{align*}
Set $y =\frac{q_2 q_3 \cdots q_{j-1}q_j q_{j+1}}{x},$ and we observe 
that for $x \in  (q_2 q_3 \cdots q_j, q_2 q_3 
\cdots q_j q_{j+1})$ we have $y \in (1, q_{j+1}).$ Thus, we obtain
\begin{align*}
h_j(y) :=   g_j \left(\frac{q_2 q_3 \cdots q_{j-1}q_j q_{j+1}}{y} \right)  &= 
\frac{q_2 q_3^{2} \cdots q_{j -1}^{j-2}
q_j^{j-1} q_{j+1}^j}{ y^{j+1}} \\
& \quad \times \bigg(1 - y + \frac{y^2}{q_{j+1}} -  
 \frac{y^3}{q_{j+1}^2q_{j}}
+ \frac{y^4}{q_{j+1}^3q_{j}^2 q_{j-1}}\bigg).  
\end{align*}
Since $q_{j+1}, q_{j}, q_{j-1} \in \left[\alpha, \frac{8}{\alpha(4 - \alpha)}\right] $ for some $\alpha \in  
 [1 + \sqrt{5}, 4),$ the polynomial in brackets satisfies the assumptions of 
Lemma~\ref {thm1:lm1}. Thus, there exists $y_j \in  (1, q_{j+1})$ such that
$h_j(y_j) \leq 0.$ Whence, there exists $x_j = \frac{ q_2 q_3 \cdots q_{j-1}q_j q_{j+1}}{y_j} \in 
(q_2 q_3 \cdots q_j, q_2 q_3 
\cdots q_j q_{j+1}),$ such that $ g_j(x_j)\leq 0.$ Taking into account (\ref{est1}),
 for every $n \geq 4$ we obtain
\begin{equation}
\label{signs2}
\forall j,  \lfloor\frac{n}{2}\rfloor +1 \leq j \leq n-1,\   \exists 
x_j  \in (q_2 q_3 \cdots q_j, q_2 q_3 
\cdots  q_{j+1})  : \     (-1)^{j-1} Q(x_j) <0.
\end{equation}

Since $q_j >1$ for all $j =2, 3, \ldots, n,$ we get $1 < q_2 < q_2q_3 <
q_2q_3q_4 < \ldots < q_2q_3 q_4 \cdots q_n, $ whence 
$x_1 < x_2 < \ldots < x_{n-1}.$ By (\ref{signs1}) and (\ref{signs2})
we have
$$Q(0) >0, - Q(x_1) >0, Q(x_2) >0, -Q(x_3) >0, \ldots ,$$
$$ (-1)^{n-1}
Q(x_{n-1})>0, (-1)^n Q(+\infty) >0.  $$
Thus, we have proved that all the zeros of $Q$ are real.

Theorem \ref{thm1} is proved.

\begin{remark}
Assumptions on $q_j$ in Theorem~\ref{thm1} could be 
slightly weakened for entire functions and polynomials of higher 
degrees if we obtain an analogue of Lemma \ref {thm1:lm1}
for polynomials of even degrees that are greater than $4.$ As it is shown 
in the proof of Lemma~\ref{thm1:lm1}, an important role in such 
considerations is played by special polynomials which have the following 
property: $\alpha =q_2 =q_4= q_6 =\ldots,$ and $\beta =q_3 =q_5= 
q_7 =\ldots,$ when $\alpha < \beta.$

The paper~\cite{NgV2021b} by T.H.~Nguyen and A.~Vishnyakova studies the 
entire functions with alternating second quotients of Taylor coefficients.
Let $f_{a, b}(x)= 1 - x +\sum_{k=2}^\infty \frac{  (-1)^k x^k}
{q_2^{k-1} q_3^{k-2} \cdots  q_k},$  be an 
entire function such that  $q_2 = q_4 = q_6 = \ldots = \alpha,$ 
$q_3 = q_5 = q_7 = \ldots = \beta,$ and $1 < \alpha < \beta.$ In \cite{NgV2021b}, it is proved that the function 
$f_{a, b}$ belongs to the Laguerre--P\'olya  class  if and only if there exists 
$x_0 \in [0,  q_2]$ such that $f_{a,b}(x_0) \leq 0.$ In addition, it is proved that if 
the function $f_{a, b}$ belongs to the Laguerre--P\'olya  class, then $ \alpha \geq q_\infty.$ 
\end{remark}

\end{document}